\documentclass[12pt,draftcls,onecolumn]{IEEEtran}
\usepackage{mathrsfs,amsmath,latexsym,amssymb}

\newtheorem{thm}{Theorem}[section]

\newcommand{\argmax}{\mathop{\mbox{\rm arg\,max}}}

\newcommand{\ps}{\mathrm{ps}}
\begin{document}
\sloppy

\title{On-line Policy Iteration with Policy Switching for Markov Decision Processes}
\author{Hyeong Soo Chang
\thanks{H.S. Chang is with the Department of Computer Science and Engineering at Sogang University, Seoul 121-742, Korea. (e-mail:hschang@sogang.ac.kr).}%
}

\maketitle
\begin{abstract}
Motivated from Bertsekas' recent study on policy iteration (PI) for solving the problems of
infinite-horizon discounted Markov decision processes (MDPs) in an on-line setting,
we develop an off-line PI integrated with a multi-policy improvement method of
policy switching and then adapt its asynchronous variant into on-line PI algorithm that
generates a sequence of policies over time. The current policy is updated into the
next policy by switching the action only at the current state while ensuring the monotonicity
of the value functions of the policies in the sequence. Depending on MDP's state-transition structure,
the sequence converges in a finite time to an optimal policy for an associated local MDP.
When MDP is communicating, the sequence converges to an optimal policy for the original MDP.
\end{abstract}


\section{Introduction}

Consider a discrete-time stochastic system whose underlying state dynamics follows a Markov decision
process (MDP). At time $k\geq 0$, the system is at a state $x_k$ in a finite set $X$ and 
an action $\pi_{k+1}(x_k)$, in the admissible
action set $A(x_k)$ of $x_k$ contained in a finite set $A$,
is taken according to the given policy $\pi_{k+1}: X\rightarrow A$ where $\pi_{k+1}(x) \in A(x)$ for all $x$ in $X$.
We denote the set of all possible policies as $\Pi$.
Once $\pi_{k+1}(x_k)$ is taken, the system makes a transition to $x_{k+1}$ by the probability 
of $P_{x_{k}x_{k+1}}^{\pi_{k+1}(x_k)}$ and obtains a reward of $R(x_{k},\pi_{k+1}(x_k))$
from the reward function $R:X\times A\rightarrow \Re$.
The transition probability is associated with the $|X|$ by $|X|$ transition-matrix $P_{\pi_{k+1}}$
whose row for $x_k$ corresponds to a probability distribution over $X$ by taking $\pi_{k+1}(x_k)$.

The problem we are concerned with is to devise an algorithm that generates a sequence of policies $\{\pi_{k+1}, k\geq 0\}$ as the system evolves from time 0 such that with an arbitrary $\pi_0\in \Pi$, the sequence of the \emph{value functions} $\{V^{\pi_{k+1}}\}$ is \emph{non-decreasing}
where for $x\in X$,
\[
   V^{\pi_{k+1}}(x) := E \biggl [\sum_{t=0}^{\infty} \gamma^{t} R(X_t,\pi_{k+1}(X_t)) | X_0 = x \biggr ],
\] where $\gamma \in (0,1)$ and $X_t$ denotes a random variable for state at $t$ by following $\pi_{k+1}$. In other words, $\{V^{\pi_{k+1}}\}$ satisfies that $V^{\pi_k} \leq V^{\pi_{k+1}}$ for all $k\geq 0$. (In what follows, any binary operator with functions or vectors over $X$ is applied componentwise.)
An important restriction is imposed on the algorithm: \emph{only}
$\pi_k(x_k)$ prescribed for $x_k$ at time $k$ is modified possibly for the next policy $\pi_{k+1}$. For all $x$ in $X\setminus \{x_k\}$, $\pi_{k+1}(x)=\pi_k(x)$. 
Because $X$ and $A$ are finite, the monotonic sequence will lead to at least a ``locally" optimal policy~\cite{bert2}. Indeed for some sufficient condition, the sequence converges to an optimal policy that achieves $\max_{\pi\in \Pi} V^{\pi}(x)$ for all $x\in X$ (see, Section~\ref{sec:on-line}).

Bertsekas recently studied an on-line policy iteration (PI) for this problem~\cite{bert1}~\cite{bert2},
where the kernel idea is the single-policy improvement equivalent to one-step application of the Newton method 
for searching roots~\cite{puterman2} (see, Section~\ref{sec:on-line} for further discussion about the rationale).
Motivated by the work,
this technical communique studies an alternative on-line PI \emph{integrated with a multi-policy improvement method} of policy switching~\cite{changps}. The algorithm is differently based on the idea of one-step application of the steepest ascent method for finding a maximum while sharing the convergence property with Bertsekas' algorithm.

\section{Off-line Policy Iteration}
\label{sec:offline}

We define \emph{switchable action set} $S^{\pi}(x)$ at $x\in X$ for $\pi\in \Pi$ as
\[
S^{\pi}(x) := \biggl \{ a \in A(x) | R(x,a) + \sum_{y\in X}P^{a}_{xy}V^{\pi}(y) > V^{\pi}(x) \biggr \}
\] and \emph{improvable-state set of} $\pi$ as $I^{\pi}  = \{ x\in X| S^{\pi}(x) \neq \emptyset \}$.
The \emph{policy improvement} result holds (see, e.g.,~\cite{puterman}):
Given $\pi \in \Pi$ with $I^{\pi} \neq \emptyset$, construct $\pi'$ such that $\pi'(x) \in S^{\pi}(x)$ for all $x\in I$ where $\emptyset \subsetneq I \subseteq I^{\pi}$ and $\pi'(x)=\pi(x)$ for all $x\in X\setminus I$. 
Then $V^{\pi'} \geq V^{\pi}$ and there exists $s\in X$ such that $V^{\pi'}(s) > V^{\pi}(s)$. 
We say that $\pi'$ \emph{strictly improves} $\pi$ and write as $\pi'> \pi$.

This result is the core of devising iterative algorithms that find an optimal policy $\pi^*$ such that $V^{\pi^*}(x) = V^*(x) := \max_{\pi \in \Pi} V^{\pi}(x)$ for all $x\in X$. Note that any $\pi\in \Pi$ that satisfies $I^{\pi} = \emptyset$ is an optimal policy.
Let $B(X)$ be the set of all possible real-valued functions defined over $X$. We introduce some operators from $B(X)$ to $B(X)$: for any $u\in B(X)$, $x\in X$, and $\pi\in \Pi$,
\[
T(u)(x) = \max_{a\in A(x)} \bigl ( R(x,a) + \gamma \sum_{y\in X} P_{xy}^{a} u(y) \bigr )
\] and
\[
T_{\pi}(u)(x) = R(x,\pi(x)) + \gamma \sum_{y\in X} P_{xy}^{\pi(x)} u(y).
\] 

The off-line PI by Howard~\cite{puterman} runs as follows: 
Set $\pi_0 \in \Pi$ arbitrarily. Loop until $I^{\pi_{n}} = \emptyset$ where $T_{\pi_{n+1}}(V^{\pi_n}) = T(V^{\pi_n}), n\geq 0$. 
Note that the terminating condition can be replaced with ``until $\pi_{n+1}=\pi_n$" if the consistent 
selection is made such that for all $x\in X$, $\pi_{n+1}(x) := \pi_n(x)$ if $x\in I^{\pi_n}$ and $\pi_n(x) \in S^{\pi}(x)$.
This algorithm, referred to as Howard-PI, is \emph{synchronous} in that \emph{all} states in $X$ are examined for improvement.

Because the strict improvement holds as long as $\emptyset \subsetneq I \subseteq I^{\pi}$,
an \emph{asynchronous} version is immediate where at every iteration a \emph{singleton} set contained
in $I^{\pi}$ is selected.
In particular, Ye~\cite{ye} studied Simplex-PI equivalent to the simplex algorithm
running with the highest-pivot rule on a linear program (LP) induced from the MDP model:
Set $\pi_0\in \Pi$ arbitrarily. Loop until $I^{\pi_{n}} = \emptyset$ with 
$\pi_{n+1}(x) = \pi_n(x)$ for all $x \in X\setminus \{s\}$, where $s$ is in
$\argmax_{x\in I^{\pi}} T(V^{\pi})(x) - V^{\pi}(x)$, and 
$T_{\pi_{n+1}}(V^{\pi_n})(s) = T(V^{\pi_n})(s), n\geq 0$. 
It can be seen that at each iteration, only
a single state (``maximally advantages state" corresponding to the highest-pivot in LP) is 
chosen for the strict policy improvement.
Simplex-PI is inappropriate to be adapted into an algorithm in our on-line setting because a specific improvable-state needs to be found at each system time, which might be different from the current state of the system, even if the current state is an improvable-state.

Suppose that an arbitrary single state, instead of a maximally advantages state, 
is chosen from $I^{\pi_n}$ and runs as in Simplex-PI. In what follows, we refer to this
asynchronous version of Howard-PI as Newton-PI. (This is intended to highlight that
the policy improvement of Howard-PI is equivalent to the 
Newton-Kantorovich iteration applied to the functional 
equation of dynamic programming for infinite-horizon discounted MDPs~\cite{puterman2} 
and switching is not involved with a particular state as in Simplex-PI.)
Bertsekas studies an algorithm~\cite{bert2} that applies Newton-PI in the
on-line setting. We will discuss about this in more detail in Section~\ref{sec:on-line}.

\subsection{Policy Switching Policy Iteration}

We now present a novel synchronous off-line PI integrated with policy switching
and derive an asynchronous version which can be readily employed via simulation 
in the on-line setting.

Let the set of \emph{all strictly better policies than} $\pi\in \Pi$ \emph{obtainable from} $I^{\pi}$ be $\beta^{\pi}$ where if $I^{\pi} = \emptyset$, $\beta^{\pi} = \emptyset$ and
\begin{eqnarray*}
\lefteqn{\beta^{\pi} = \bigl \{\pi'\in \Pi | \forall I \in 2^{I^{\pi}}\setminus \{\emptyset\},
  \forall x\in I \mbox{ } \pi'(x) \in S^{\pi}(x)} \\ 
& & \hspace{3cm} \mbox{ and } \forall x\in X\setminus I \mbox{ } \pi'(x)=\pi(x) \bigr \}
\end{eqnarray*} otherwise. In other words, for non-optimal $\pi$, any $\pi'\in \beta^{\pi}$ strictly improves $\pi$.
Obviously each policy obtained after applying the policy improvement step of Howard-PI, Simplex-PI, and Newton-PI to non-optimal $\pi$, respectively, is in $\beta^{\pi}$. 

Given a nonempty $\Delta\subseteq \Pi$, define \emph{policy switching} with $\Delta$ as $\pi_{\ps}(\Delta) \in \Pi$ such that for all $x\in X$,
\[ \pi_{\ps}(\Delta)(x) = \phi_x(x) \mbox{ where } \phi_x\in \argmax_{\pi\in\Delta} V^{\pi}(x).
\] Chang \emph{et al.}~\cite{changps} showed $V^{\pi_{\ps}(\Delta)}\geq V^{\phi}$ for all $\phi \in \Delta$.
The \emph{multi-policy improvement} result immediately follows: 
\begin{thm}
\label{thm1} 
For any nonempty $\Delta\subseteq \Pi$ and $\pi\in \Pi$, if $\Delta \cap \beta^{\pi} \neq \emptyset$, then $\pi_{\ps}(\Delta) > \pi$ and $V^{\pi_{ps}(\Delta)}\geq V^{\phi}$ for all $\phi \in \Delta$.
\end{thm} We remark that a related result exists based on ``parallel rollout" and 
a generalized Howard-PI by using the result was studied in~\cite{changpsi}, 
where at each iteration the multi-policy improvement is applied over a given set
that includes the policy generated from a previous iteration. However, it seems difficult to 
transform the algorithm into an on-line version because preserving the multi-policy improvement 
of parallel rollout in the on-line setting appears problematic. In fact, this obstacle lead to
one of the motivations for considering the policy-switching approach here.

We provide a generic synchronous version of policy-switching PI (PSPI):
Set $\pi_0\in \Pi$ and $\Delta_0\subseteq \Pi$ arbitrarily such that $\Delta_0 \cap \beta^{\pi_0}\neq \emptyset$.
Loop until $I^{\pi_{n}} = \emptyset$ with $\pi_{n+1} = \pi_{\ps}(\Delta_n)$. Set $\Delta_{n+1}$ such that $\Delta_{n+1} \cap \beta^{\pi_{n+1}}\neq \emptyset, n\geq 0$.
The convergence to an optimal policy is obvious because $\pi_{n+1} > \pi_n, n\geq 0$.
The generality of $\{\Delta_n\}$ provides broad design-flexibility and in particular $\{\Delta_n\}$ can be set to include the policies in the sequences generated by other PI algorithms run in parallel (with the same initial policy), in which case PSPI converges no slower than any other PI algorithms in terms of iterations.

We now move to an asynchronous version of PSPI. For $\pi\in \Pi$ such that $I^{\pi}\neq \emptyset$ and $x\in I^{\pi}$, let
\[
   \beta^{\pi}_x = \bigl \{\pi'\in \Pi | \forall s\in X\setminus\{x\}, \pi'(s) = \pi(s) \bigr \}.
\] Trivially, there are exactly $|A(x)|$ policies of $\pi_a, a\in A(x),$ in $\beta^{\pi}_{x}$ such that for each $a\in A(x)$, $\pi_a(x)=a$. 
Note that $\pi$ itself is in $\beta^{\pi}_x$ and there exists $\pi'\in \beta^{\pi}_x$ that satisfies $T_{\pi'}(V^{\pi})(x) = T(V^{\pi})(x)$. It follows that
because
\[
   V^{\pi_{\ps}(\beta^\pi_x)} \geq V^{\phi}, \forall \phi \in \beta^{\pi}_x,
\] the policy switching policy with $\beta^\pi_x$ is no worse than $\pi'$ generated by Newton-PI.

At this point, a careful attention needs to be paid to the classic result that states that $\pi'$ defined as
$T_{\pi'}(V^{\pi})=T(V^{\pi})$ in Howard-PI strictly improves non-optimal $\pi$.
It \emph{does not} say that $\pi'$ improves any other policy that strictly improves $\pi$.
In $\beta^\pi_x$, there can be more than one policies that strictly improve $\pi$
because $|S^{\pi}(x)|>1$ in general. Even there can exist $\phi\in \beta^\pi_x$ that 
satisfies $T_{\phi}(V^{\pi}) = T(V^{\pi})$ but $\pi'\neq \phi$.
It is the \emph{main profit} of policy switching that $\pi_{\ps}(\beta^\pi_x)$ is no worse than any policy that strictly improves $\pi$ in $\beta^{\pi}_x$.

Suppose that $\pi_{a'}$ in $\beta^{\pi}_x$ for some $a'\in A(x)$ satisfies that $T_{\pi_{a'}}(V^{\pi})\leq V^{\pi}$ and $\pi(x)\neq a'$. Then $V^{\pi_{a'}}\leq V^{\pi}$. Therefore, we can see that if $\pi(x)$ is consistently selected as $\pi_{\ps}(\Delta)(x)$ whenever $\pi \in \argmax_{\phi \in \Delta} V^{\phi}(x)$,
such non-improving policy $\pi_{a'}$ of $\pi$ is always excluded in the switching process.

The above discussion naturally induces an asynchronous PSPI
on which the on-line PSPI in Section~\ref{sec:on-line} is based on:
Set $\pi_0\in \Pi$ arbitrarily. Loop until $I^{\pi_{n}} = \emptyset$ with 
$\pi_{n+1}(x_n) = \pi_{\ps}(\beta^{\pi_n}_{x_n})(x_n)$ for arbitrarily chosen state $x_n\in I^{\pi_n}$ 
and $\pi_{n+1}(x) = \pi_n(x)$ for all $x\in X\setminus \{x_n\}, n\geq 0$.
The convergence to an optimal policy is again guaranteed due to the consecutive strict
improvements. 

\section{On-line Policy Iteration with Policy Switching}
\label{sec:on-line}

Bertsekas' main motivation for developing on-line PI is from the supportive success of the AlphaZero program (an improvement of the breakthrough program AlphaGo~\cite{changgo}) 
that plays Go, Chess, and some other games (see, e.g.,~\cite{silver} and~\cite{bert1} and the references therein).
AlphaZero's self-learning (without any help from (human) experts) in its off-line training phase is basically applying a variant of off-line PI and when AlphaZero
plays on-line against another player, it applies a policy improvement through its on-line play
due to inherent ``approximations" during the off-line training~\cite{bert1}.
The underpinning critical idea of Bertsekas is that because
Howard-PI is equivalent to applying the Newton method within MDP setting, an application of the policy improvement for the next policy (corresponding to drawing the tangent line to the curve of interest at the current point on the curve and finding the crossing point with the axis associated with roots of the curve)
would 
subsume the well-known advantage of the Newton method:
If the starting point is selected ``well", then just a single iteration provides a good approximation to the solution (or zero).
This advocates in particular his simulation-based approach for the policy improvement by \emph{rolling out} (simulating) the current (base) policy in the MDP systems (over a selected finite horizon for approximation of infinite horizon).
There are already some works in the literature that report positive results of using the rollout approach (see, e.g.,~\cite{bert1} and the references therein).

However, we argue that some limitations exist for the rollout approach.
There seems \emph{no general rule for selecting an initial point} for the Newton method. 
Some theories characterize the convergence (rate) of the Newton method with some
property on the initial point but do not summarize into a general rule and more importantly,
the conditions are not directly applicable and verifiable within the MDP context 
(see, e.g.,~\cite{galanti}).
Because $V^*$ is not known in advance in general, how well a current policy approximates $V^*$ is not known either.
Indeed, there is no known general rule for selecting $\pi_0$ while applying Howard-PI.
Even if the upper bounds on the time-complexity of Howard-PI in terms of iterations have been studied (see, e.g.,~\cite{scher} and the references therein), there is no known upper bound yet that characterizes the structure of the initial policy chosen in the bound.

All of existing results in the literature (including AlphaZero) about the rollout approach are mainly \emph{experimental} to the author's best knowledge. 
Consider two policies $\pi_1$ and $\pi_2$ such that
$T_{\pi_1}(V^{\phi_1})=T(V^{\phi_1})$ and $T_{\pi_2}(V^{\phi_2})=T(V^{\phi_2})$. Even if $V^{\phi_1}\geq V^{\phi_2}$, there seems no known result that allows to say that $V^{\pi_1}\geq V^{\pi_2}$. In other words, \emph{even if we roll out a better base-policy} $\phi_1$ than the current one $\phi_2$,
it is difficult or even impossible to say that the resulting policy $\pi_1$ from rolling out $\phi_1$ is better than the policy $\pi_2$ from rolling out $\phi_2$.
\emph{What can be guaranteed} is that if we roll out a non-optimal $\pi_0\in \Pi$ as a base policy, the resulting policy $\pi_1$ strictly improves \emph{only} $\pi_0$. Necessarily, we then should roll out $\pi_1$ in order to be sure that we improve not only $\pi_1$ but also $\pi_0$.
A \emph{consecutive on-line application of the policy improvement} is required over time to ensure the monotonicity over all past policies.

To implement the rollout approach, we need to approximate \emph{double expectation} besides maximization or minimization over the action set. One expectation is for the next-state distribution and another is for the value function of the base policy to be rolled out.
This is rather cumbersome because we need to sample next states first and for each sampled next state, we need to rollout the base policy starting from the next state.

Still, the rollout approach is expected to be a viable and good starting-solution methodology when we face with the curse of the dimensions in the model but we have some heuristic(s) or expert knowledge, etc. to be rolled out.
Formally, Bertsekas-OPI (on-line PI by Bertsekas) generates a sequence $\{\pi_{k+1}, k\geq 0\}$ of policies with arbitrarily chosen $\pi_0$ 
by applying the policy improvement of Howard-PI consecutively over time only at the current state $x_k$ of the system such that $T_{\pi_{k+1}}(V^{\pi_k})(x_k) = T(V^{\pi_k})(x_k)$ and 
set $\pi_{k+1}(x)=\pi_k(x), x\in X\setminus\{x_k\}$. It is simply an asynchronous Newton-PI over $\{x_k\}$. The system makes a transition to a random next state $x_{k+1}$ according to $P_{\pi_{k+1}}$ by taking $\pi_{k+1}(x_k)$.
Because $x_k$ at some $k$ can be non-improvable state and a probabilistic transition is made by $\pi_{k+1}(x_k)$, the convergence does not necessarily hold as before.

Recall that we consider MDPs with a finite state set and a finite action set.
Bertsekas~\cite{bert2} established for the finite MDPs that for $\{\pi_k\}$ by Bertsekas-OPI,
there exists $k'<\infty$ such that $\pi_k = \pi$ for all $k>k'$.
In particular, $\pi$ is an optimal policy for the local MDP problem induced from
the original MDP with the same $X$, $R$, and $P_{\phi}, \phi\in \Pi$, but $A$ is
replaced by $A_{\chi}$ for some $\chi\subseteq X$, where
for all $x\in \chi$ and all $y\in X\setminus \chi$, $P^{\pi(x)}_{xy}=0$, with 
$A_{\chi}(x) = A(x)$ for $x\in \chi$ and for $x\in X\setminus \chi$, $A_{\chi}(x) = \pi(x)$. 
In other words, $\{\pi_k\}$ converges in a finite time to 
a locally optimal policy $\pi$ by visiting the states in $\chi$ infinitely often within $\{x_k\}$ 
such that $\pi$ solves the local MDP over $\chi$ restricted by taking $\pi(x)$ when $x$ is not in $\chi$.
This convergence property brings up an immediate sufficient condition:
When MDP is \emph{communicating}~\cite{kallen}, the convergence in a finite time to an optimal policy 
for the original MDP is achieved.

We are now back to PSPI.
On-line PSPI is adapted directly from the asynchronous PSPI:
Starting with arbitrary $\pi_0\in \Pi$, generate $\{\pi_{k+1},k\geq 0\}$ such that
$\pi_{k+1}(x_k) = \pi_{\ps}(\beta^{\pi_k}_{x_k})(x_k)$
and $\pi_{k+1}(x) = \pi_k(x)$ for all $x\in X\setminus \{x_k\}$ and take $\pi_{k+1}(x_k)$ at $x_k$.
Because $\pi_k \in \beta^{\pi_k}_{x_k}, k\geq 0$, the following monotonicity result
that makes on-line PSPI converge in a finite time to an optimal policy as in Bertsekas-OPI holds:
\begin{thm}
For $\{\pi_{k+1},k\geq 0\}$ generated by on-line PSPI, $V^{\pi_{k+1}} \geq V^{\pi_k}$ for all $k\geq 0$.
\end{thm} 

The on-line PSPI inherits the spirit of the rollout approach but rolling out \emph{multiple} policies.
Starting with $x_k$, we roll out each $\pi_a, a\in A(X),$ (independently in parallel) in order to estimate $V^{\pi_a}(x_k)$ (by using the same random number sequences if possible). 
Note that each $\pi_a$ simulation corresponds to estimating \emph{only one} expectation unlike Bertsekas-OPI. Simply take the action of a policy that maximizes the estimates obtaining $\pi_{\ps}(\beta^{\pi_k}_{x_k})(x_k)$.
Note that the probability that this simple Monte-Carlo method, which is basically equivalent to the method of sample average approximation (SAA)~\cite{kleywegt}, finds a maximizing policy approaches one \emph{exponentially fast} in the number of the sample paths of the rollout.
The problem of obtaining $\argmax_{\pi_a, a\in A(x_k)}V^{\pi_a}(x_k)$ is also equivalent to the best arm identification~\cite{bubeck} in stochastic multi-armed bandits (MABs) by viewing $\pi_a$ as an arm and
the value function value at the current state $V^{\pi_a}(x_k)$ as the expected performance of the arm.
In this context, on-line PSPI is faced with a sequence of MAB problems and available MAB algorithms for finding a best arm can be also applied at each time (see, e.g.,~\cite{bubeck}).
Note that these reductions to SAA and MAB problems are possible due to the methodology of estimating one expectation unlike double expectation in the rollout approach.

The idea behind policy switching can be attributed to approximating the steepest ascent direction 
while applying the steepest ascent algorithm. At the current location $\pi_k$,
we find ascent directions relative to $\pi_k$ over the local neighborhood of $\beta^{\pi_k}_{x_k}$
by $V^{\phi}(x_k)-V^{\pi_k}(x_k), \phi\in \beta^{\pi_k}_{x_k}$.
A steepest ascent direction, $\pi_{\ps}(\beta^{\pi_k}_{x_k})$, is then obtained by ``combining" all of the possible ascent directions. In particular, the ascent direction $\phi$ that satisfies $T_{\phi}(V^{\pi_k})(x_k)=T(V^{\pi_k})(x_k)$ is included while combining.

\section{Concluding Remark}

Comparing on-line PSPI with Bertsekas-OPI in terms of the convergence rate is nontrivial because analyzing the relative distances to $V^*$ from the policies in the paths generated by both algorithms are difficult and the state sequences are different in general even if the initial policies are the same.
PSPI here should be regarded as an alternative method to Bertsekas-OPI but not competitive one.

If desired, $\beta^{\pi_k}_{x_k}$ in PSPI can be replaced with some $\Phi_k = \beta^{\pi_k}_{x_k} \cup \Delta_k$, where $\Delta_k$ can include any policy in $\Pi$ such as (human) experts knowledge or heuristics, if represented as policies, and even some policy by running Bertsekas-OPI in parallel, etc. We argue that PSPI has a significant extensibility or flexibility in designing on-line algorithms for solving MDPs on-line.

\end{document}